\documentclass[12pt]{amsart}
\let\oldlabel=\label
\def\prellabel
   {\marginparsep=1em\marginparwidth=44pt
\def\label##1{\oldlabel{##1}
\ifmmode
\else
  \ifinner
  \else \marginpar{{\footnotesize\ \\\normalshape\tt ##1}}
  \fi
\fi}}
\def\Q{\ifhmode\textqed\fi
   \ifmmode\ifinner\quad\Qsymbol\else\dispqed\fi\fi}
\def\textqed{\unskip\nobreak\penalty50
    \hskip2em\hbox{}\nobreak\hfil\Qsymbol
    \parfillskip=0pt \finalhyphendemerits=0}
\def\dispqed{\rlap{\qquad\Qsymbol}}

\usepackage{amscd}

\def\mm{{\mathfrak m}}  
\def\nn{{\mathfrak n}}

\def\and{\operatorname{and}}

\def\max{\operatorname{max}}
\def\Hom{\operatorname{Hom}}

\def\mm{{\frak m}}

\def\bc{\begin{corollary}}
\def\ec{\end{corollary}}
\def\ben{\begin{enumerate}}
\def\een{\end{enumerate}}
\def\beqn{\begin{eqnarray*}}
\def\eeqn{\end{eqnarray*}}
\def\bd{\begin{definition}}
\def\ed{\end{definition}}
\def\bl{\begin{lemma}}
\def\el{\end{lemma}}
\def\bpf{\begin{proof}}
\def\epf{\end{proof}}
\def\bp{\begin{proposition}}
\def\ep{\end{proposition}}
\def\bt{\begin{theorem}}
\def\et{\end{theorem}}
\def\bex{\begin{example}}
\def\eex{\end{example}}

\let\epsilon=\varepsilon
\let\phi=\varphi
\let\theta=\vartheta

\newtheorem{lemma}{Lemma}[section]
\newtheorem{corollary}[lemma]{Corollary}
\newtheorem{theorem}[lemma]{Theorem}
\newtheorem{proposition}[lemma]{Proposition}
\theoremstyle{definition}
\newtheorem{definition}[lemma]{Definition}

\newtheorem{example}[lemma]{Example}

\advance\headheight1.15pt

\begin{document}
\title {On a Theorem of Macaulay on Colons of Ideals$^*$}
\thanks{\small $^*$Presented at the Third National Conference 
on Commutative Algebra and Algebraic Geometry, Indian Institute of
Science, Bangalore, India, October 16-21, 2000}

\author{ J. K. Verma}
\address{Department of Mathematics, IIT Bombay, Powai, Mumbai, India -
400076}
\email{jkv@math.iitb.ac.in}

\maketitle

\begin{center}
{\em \large Dedicated to my teacher William Heinzer}
\end{center} 

\begin{abstract}

Let $R$ be a polynomial ring over a field of dimension $n.$ Let $\mm$ be
the maximal
homogeneous ideal of $R.$ Let   $I$ be a
complete intersection homogeneous ideal of $R$ minimally generated by
polynomials
of degrees $ d_1, d_2, \ldots, d_n.$
F. S. Macaulay   showed that for all integers $i=0,1,\ldots, \delta+1,$
$I:\mm^i=I+\mm^{\delta+1-i}$ where  $\delta= d_1+d_2+ \ldots + d_n-n.$ 
We provide a modern proof of a  generalization of this  result to
standard graded Gorenstein rings. 

\end{abstract}

\maketitle
 
\setcounter{section}{1}
In \S 86 of his monograph \cite{mac} of 1916, F. S. Macaulay
proved the following

\bt

Let $R$ be the polynomial ring $k[X_1, X_2, \ldots, X_n]$ over a field
$k$. Let $I$ be  a height $n$ ideal  of $R$ generated by homogeneous
polynomials
$f_1, f_2, \ldots, f_n$ of degrees 
$ d_1, d_2, \ldots, d_n$ respectively.
Set $\mm=(X_1, X_2, \ldots, X_n)R$ and 
$\delta= d_1+d_2+ \ldots + d_n-n.$ 
Then 
for all integers $i=0,1,\ldots, \delta+1,$

$$
I:\mm^i=I+\mm^{\delta+1-i}.
$$
\et

 In his expository paper \cite{gri}  P. Griffiths  used
Macaulay's theorem
in Hodge theory of smooth hypersurfaces in projective space. He also
provided a proof of this theorem  by using the local duality theorem in
the complex case involving the Grothendieck residue symbol. The objective
of this note is to provide a modern algebraic proof which seems to be
lacking in the literature. We will  generalize Macaulay's theorem to 
$\mm$-primary ideals in standard graded  Gorenstein  rings.

By taking the quotient of $R$ by $I$, it is enough to prove the theorem
for standard graded Artinian Gorenstein rings.

A crucial result used  in our proof is a theorem of Macaulay about Hilbert
series of Gorenstein graded algebras. For a modern proof, see 
Corollary
4.4.6 of \cite{bh}. This has been generalized to Gorenstein graded rings
over Artinian rings in \cite{jv}. We shall use  the zero dimensional case
of this. We recall the statement for the convenience of the reader. Let
$\ell(.)$ denote length. 

\bt
Suppose that $S=\bigoplus_{n=0}^{\infty}S_n$ is a standard
graded Artinian Gorenstein algebra over an Artinian local ring $S_0.$ Let 
$\delta = \max\{ n \;|\; S_n \neq 0\}.$ Then 
$\ell(S_i)=\ell(S_{\delta-i})$ for all $i=0, 1, \ldots, \delta.$ 
\et

We  recall certain standard facts from \cite{bh} about local cohomology.
Let $R=\oplus_{n=0}^{\infty}R_n$ be a $d$-dimensional finitely generated
standard graded
algebra over a Noetherian local ring $R_0.$
 Let $\mm$ denote the irrelevant ideal
$\oplus_{n=1}^{\infty}R_n$. The $a$-invariant of $R$ is defined by

$$a(R)=\max \{ n | \left[H^{d}_{\mm} (R)\right]_n \neq 0\}.$$
If $R$ is Cohen-Macaulay and  $f$ is a homogeneous  nonzerodivisor
of degree $d$ then $a(R/fR)=a(R)+d.$ 
If $R$ is a polynomial ring over a field in $n$ variables then
 $a(R)=-n.$
Thus If $ {\bf f}=(f_1, f_2, \ldots , f_n)$ is a regular sequence
in a polynomial ring with
deg$(f_i)=d_i$ for all $i=1, 2, \ldots, n$ then $a(R/{\bf f})=
d_1+d_2+\ldots+d_n-n.$ This explains the value of $\delta$ in Macaulay's
theorem.

\bt
\label{main}
Let $(R,\mm)$ be an Artinian Gorenstein local ring. Let $I$ be an
$\mm$-primary ideal. Put $H(I,n)=\ell(I^n/I^{n+1}).$ Let 
$\delta=\max\{i | I^i \neq 0\}.$  Then $0:I^i=I^{\delta +1-i}$
for all  $i=0,1,\ldots, \delta$ if and only if $H(I,i)=H(I,\delta-i) $
for all  $i=0,1,\ldots, \delta$.
\et

\bpf
Since $R$ is Gorenstein, it is the injective hull of its residue field.
Hence by
Matlis duality $\ell(\Hom(R/I^i,R))=\ell(R/I^i)$ for
all $i.$ Since
$I^{\delta+1-i} \subseteq (0:I^i)$ and $\ell(\Hom(R/I^i,R))=\ell(0:I^i),$
we conclude that $0:I^i=I^{\delta +1-i}$
for all  $i=0,1,\ldots, \delta$ if and only if
$\ell(R/I^i)=\ell(I^{\delta+1-i})$ for all  $i=0,1,\ldots, \delta$ 
if and only if $H(I,0)+ \ldots + H(I,i)=H(I,\delta-i)+\ldots +
H(I,\delta)$ for all  $i=0, 1, \ldots, \delta$. It is easy to see that the
last statement is equivalent to
the condition $H(I, i)=H(I,\delta-i)$ for all $i=0, 1, \ldots, \delta.$ 

 \epf

\bc 
Let $R$ be a standard graded Artinian Gorenstein ring over an Artinian
local ring $(S,\nn).$
Let $\mm $ be the graded ideal of $R$ generated by elements of positive
degree  with $\mm^{\delta} \neq 0$ and
$\mm^{\delta+1}=0.$   Then for all $i=0,1,\ldots, \delta,$
$$
0:\mm^i=\mm^{\delta+1-i}.$$
\ec  

\bpf 
Since $R$ is standard, 
the associated graded ring of $R$ with respect to $\mm$ is $R$ itself.
Since $R$ is Gorenstein,  its Hilbert function is symmetric
by \cite{jv}.

\epf

\bex
 We present an example due to U. Storch which shows that
Macaulay's theorem is not valid for Gorenstein local rings. Let $k$ be a
field of characteristic $2.$ Let $X$ and $Y$ be indeterminates. Put
$R=k[X,Y]/I$ where $I=(X^2+Y^2, X^2+XY+Y^3).$ Then $R$ is a
Gorenstein local ring. The Hilbert series
$H(G,\lambda)=\sum_{i=0}^{\infty}\ell(\mm^i/\mm^{i+1})\lambda^i $
 of the associated graded $G$ 
 of $R$ with respect to its maximal ideal $\mm$ is given by
 
$$ H(G,\lambda)=1+2\lambda+\lambda^2+\lambda^3.
$$
Thus the Hilbert function is not symmetric. Hence Macaulay's theorem
does not hold for $R$ by Theorem \ref{main}. 
\eex

\noindent
{\bf Acknowledgement :} The author thanks U. Storch and B. Singh for
useful conversations.

\end{document}